\documentclass[12pt]{iopart}
\usepackage{iopams}
\usepackage{amsthm}
\usepackage{color}
\usepackage{graphicx}
\newcommand{\Area}{\mathrm{Area}}
\newcommand{\Length}{\mathrm{Length}}
\newcommand{\Inter}{\mathrm{Int}}
\newcommand{\text}[1]{\mbox{#1}}
\newcommand{\mod}{\,\mathrm{mod}\,}
\newcommand{\pard}{{\partial_D}}
\newcommand{\parn}{{\partial_N}}
\eqnobysec
\theoremstyle{plain}
\newtheorem{thm}[equation]{Theorem}

\theoremstyle{definition}

\theoremstyle{remark}

\begin{document}

\title{Isospectral domains with mixed boundary conditions}

\author{Michael Levitin$^1$, Leonid Parnovski$^2$, Iosif Polterovich$^3$}

\address{$^1$ Maxwell Institute, Department of Mathematics, Heriot-Watt University, Riccarton, Edinburgh EH14 4AS, U.K.;
e-mail {\sffamily M.Levitin@ma.hw.ac.uk}}
\address{$^2$ Department of Mathematics, University College London, Gower Street, London WC1E 6BT, U.K.;
e-mail {\sffamily leonid@maths.ucl.ac.uk}}
\address{$^3$ D\'ept. de math\'ematiques et de statistique, Universit\'e de Montr\'eal, CP 6128 succ Centre-Ville,
Montreal, QC H3C 3J7, Canada;
e-mail {\sffamily iossif@dms.umontreal.ca}}
%
%
%

%
%
\begin{abstract} We construct a series of examples of planar isospectral domains
with mixed Dirichlet-Neumann boundary conditions. This is a modification of a classical
problem proposed by M. Kac. 


{\ }

\noindent Published in \emph{J. Phys. A: Math. Gen.} {\bf 39} 2073--2082 (2006). This version corrects the statement of Theorem 5.1.
\end{abstract}
%
\section{Introduction}\label{sec:intro}
Let us recall the classical question of Mark Kac, ``\emph {Can one
hear the shape of a drum?}" \cite{K}. He asked whether there exist two non-isometric domains 
on the plane such that the spectra of the (Dirichlet) Laplacian on them coincide (such domains are called \emph{isospectral}).
For arbitrary planar domains it was
answered negatively in \cite{GWW} using an algebraic construction
of \cite{Sun}. See also reviews and extensions \cite{BCDS, Bro,
Bus, Ch} and references therein, as well as illustrative 
numerics in \cite{BetTr}. Also, it was observed in \cite{SlHu} that the constructions of \cite{BCDS,Ch} 
work in the case of domains with fractal boundaries. 
However, to the best of our knowledge there are still no known examples
of sets of more than two non-isometric isospectral domains, as well as of non-simply connected domains.
Also, Kac's question still remains open for smooth (as well as for convex) domains. 

In the present paper we suggest the following
modification of isospectrality question for the case of mixed Dirichlet-Neumann 
boundary conditions (the so called Zaremba problem \cite{Za}) which develops the approach suggested in \cite{JLNP}.
Some related numerically constructed (though more complicated) examples can be found in \cite{DG}.

Let $\Omega_j\subset\mathbb{R}^2$, $j=1,2$ be two bounded domains, their piecewise smooth boundaries
being decomposed as $\partial\Omega_j =
\overline{\pard\Omega_j\cup\parn\Omega_j}$, where
$\pard \Omega_j$, $\parn\Omega_j$ are finite unions of open
segments of $\partial \Omega_j$ and $\pard \Omega_j \cap
\parn\Omega_j = \emptyset$. Suppose that there are no
isometries of $\mathbb{R}^2$ mapping $\Omega_1$ onto $\Omega_2$ in such a way that 
$\pard\Omega_1$ maps onto $\pard\Omega_2$. (We shall call such pairs of domains
\emph{nontrivial}.) Consider on each  $\Omega_j$ a mixed boundary value problem for the Laplacian, 
\[
-\Delta u=\lambda u\quad \text{in} \,\,\, \Omega_j\,,\quad u|_{\pard\Omega_j}=0\,,
\quad \partial u/\partial n|_{\parn\Omega_j}=0\,,
\]
and denote its spectrum by $\sigma_{DN}(\Omega_j)$.
Our aim is to study nontrivial isospectral pairs $\Omega_1$, $\Omega_2$ 
(i.e. such that $\sigma_{DN}(\Omega_1)=\sigma_{DN}(\Omega_2)$). We present a series of examples of 
such pairs and provide some necessary  conditions for
mixed isospectrality. 

\section{Basic example}\label{sec:basex} 
\begin{figure}[!hbt]
\begin{center}
\framebox[0.8\textwidth]{\includegraphics[width=0.60\textwidth]{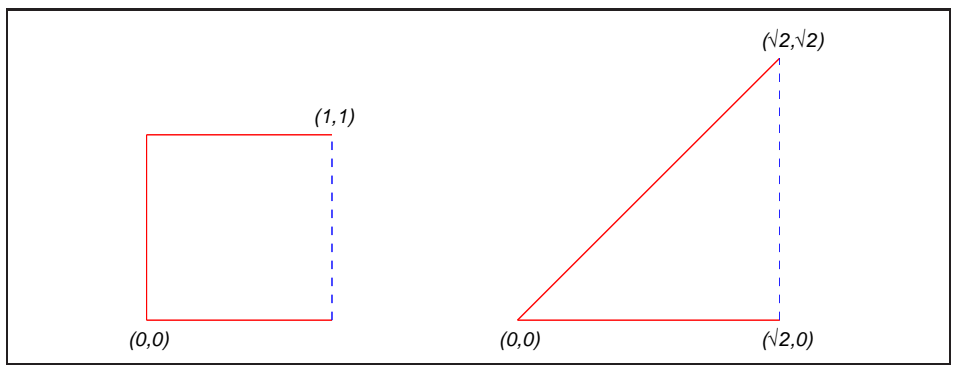}}
\caption{The unit square $\Omega_1$ and the isosceles triangle $\Omega_2$. Here and further on, 
\textcolor{red}{red solid lines \full} denote the Dirichlet boundary
conditions and \textcolor{blue}{blue dashed lines  \dashed} denote the Neumann
ones.\label{fig:trivialex}}
\end{center}
\end{figure}

Let $\Omega_1=(0,1)^2$ be a unit square with 
$\parn{\Omega_1}=\{1\}\times(0,1)$ and $\pard{\Omega_1}= \{0\}\times[0,1]\cup
(0,1)\times\{0,1\}$. Let $\Omega_2$ be an isosceles right-angled triangle 
$\{(x,y)\in\mathbb{R}^2: 0<x<\sqrt{2}\,,\ 0<y<x\}$ with 
$\parn{\Omega_2}=\{\sqrt{2}\}\times(0,\sqrt{2})$ and 
$\pard{\Omega_2}=\{0\}\times(0,\sqrt{2})\cup\{(x,x): x\in(0,\sqrt{2})\}$ (see
Figure~\ref{fig:trivialex})

The spectra $\sigma_{DN}(\Omega_j)$ are easily calculated by separation of variables. 
The eigenfunctions for $\Omega_1$ are
\[
\sin((1/2+m)\pi x)\sin(n\pi y)\,,\qquad\text{for }n=1,2,\dots\,,\ m=0,1,2,\dots\,,
\]
and the eigenfunctions for $\Omega_2$ are
\begin{eqnarray}
\fl\sin\left(\frac{(1/2+k)\pi x}{\sqrt{2}}\right)\sin\left(\frac{(1/2+l)\pi y}{\sqrt{2}}\right)&-
\sin\left(\frac{(1/2+l)\pi x}{\sqrt{2}}\right)\sin\left(\frac{(1/2+k)\pi y}{\sqrt{2}}\right)\,,\nonumber\\&
\text{for }k=0,1,2,\dots\,,\ l=0,1,2,\dots\,,\ k>l.
\end{eqnarray}
The corresponding spectra (with account of multiplicities) are 
$\sigma_{DN}(\Omega_1)=\{\lambda_{m,n}\}$ and $\sigma_{DN}(\Omega_2)=\{\mu_{k,l}\}$ with
\begin{equation}
\lambda_{m,n}=\frac{\pi^2}{4}\,\left((2m+1)^2+4n^2\right)
\end{equation} 
and
\begin{equation}
\mu_{k,l}=\frac{\pi^2}{4}\,\frac{(2k+1)^2+(2l+1)^2}{2}\,.
\end{equation} 
\begin{thm}
\[
\sigma_{DN}(\Omega_1) =\sigma_{DN}(\Omega_2)
\]
\end{thm}

\begin{proof}  Indeed, it is easy to check that 
\[
\mu_{k,l}=
\cases{
\lambda_{j,l+j+1}&if $k-l=2j+1$\,,\\
\lambda_{l+j,j}&if $k-l=2j$\,.
}
\]
On the other hand,
\[
\lambda_{m,n}=
\cases{
\mu_{m+n,m-n}&if $m\ge n$\,,\\
\mu_{m+n,n-m-1}&if $m<n$\,.
}
\]
These two correspondences establish a bijection between $\sigma_{DN}(\Omega_1)$ and $\sigma_{DN}(\Omega_2)$.
\end{proof}

The example above shows that isospectral domains with mixed boundary conditions can be quite simple 
(compared to classical Dirichlet isospectral pairs in \cite{GWW}, \cite{BCDS} and numerical examples in \cite{DG}).
Indeed, dependence of the spectra on  boundary decomposition brings more flexibility 
to  the problem. See also section~\ref{sec:more} for other illustrations of this phenomenon. 

We note that this example is somewhat reminiscent of Chapman's example \cite{Ch} of two disconnected Dirichlet 
isospectral domains: in his case the first disconnected domain is a disjoint union of a square of side one and an isosceles right 
triangle of side two, and the second one is a disjoint union of a rectangle with sides one and two and an isosceles right triangle of 
side $\sqrt{2}$.

\section{Main construction}\label{sec:main} 
Example of section~\ref{sec:basex} is in fact the easiest implementation of the following algorithm 
for the construction of isospectral domains which we outline below. 

We start by describing a suitable class of ``construction blocks'' which we shall later use to build  pairs of
planar isospectral domains. Let $a$, $b$ be two lines on the plane (which may be parallel), and  let $K$ be a
bounded open set lying in a sector formed by $a$ and $b$ (or between them if they are parallel). $ K$ need not
be connected, but we assume that $\partial K$ has non-empty intersections with $a$ and $b$  which we denote
$\partial_a K$ and $\partial_b K$, respectively. Let  $\partial_0 K:=\partial K\setminus(\partial_a K\cup
\partial_b K)$ be the remaining part of the boundary $\partial K$.

Let $T_{a}$, $T_b$ denote the reflections with respect to the lines $a$ and $b$. We first construct the 
domains $\Omega_1$ and $\Omega_2$ just by adding to $K$ its image under reflections $T_a$, $T_b$, 
respectively:
\[
\Omega_1:=\Inter(\overline{K\cup T_a(K)})\,,\qquad
\Omega_2:=\Inter(\overline{K\cup T_b(K)})\,.
\]

\begin{figure}[!hbt]
\begin{center}
\framebox[0.8\textwidth]{\includegraphics[width=0.60\textwidth]{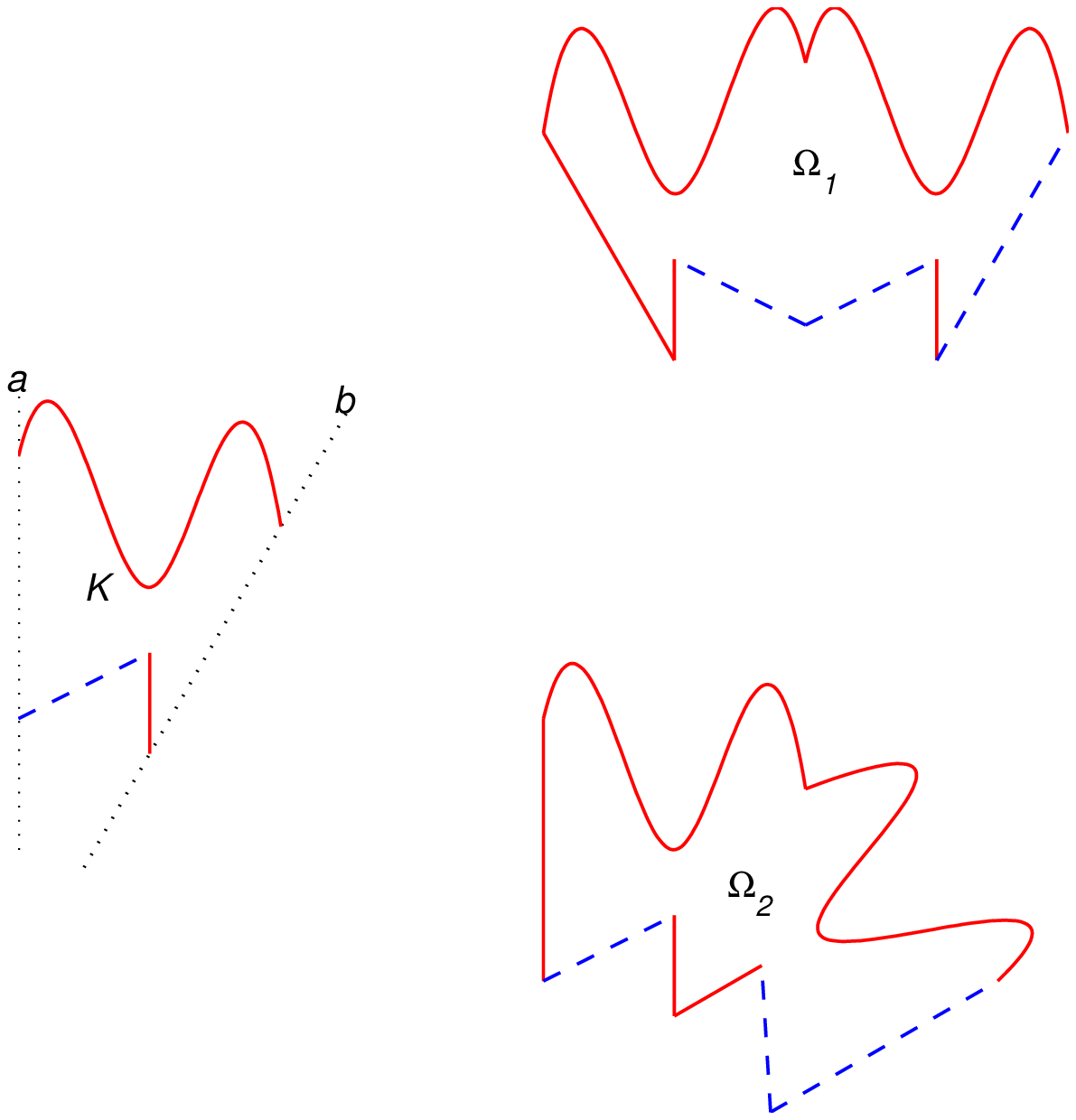}}
\caption{A generic construction block and resulting domains $\Omega_1$ and
$\Omega_2$.\label{fig:constblock}}
\end{center}
\end{figure}

We now need to impose mixed boundary conditions on $\Omega_1$ and 
$\Omega_2$. To do so, let us first decompose $\partial_0 K$ into the union of two non-intersecting 
sets $\partial_{0,D} K$ and $\partial_{0,N} K$ (one of which may be empty).
Then we set (see figure~\ref{fig:constblock})
\[
\partial_D \Omega_1:=\partial_{0,D} K \cup T_a(\partial_{0,D} K) \cup T_a(\partial_b K)\,;
\]
\[
\partial_N \Omega_1:=\partial_{0,N} K \cup T_a(\partial_{0,N} K) \cup \partial_b K\,;
\]
\[
\partial_D \Omega_2:=\partial_{0,D} K \cup T_b(\partial_{0,D} K) \cup \partial_a K\,;
\]
\[
\partial_N \Omega_2:=\partial_{0,N} K \cup T_b(\partial_{0,N} K) \cup T_b(\partial_a K)\,;
\]

\begin{thm}\label{thm:main}
For any choice of lines $a$, $b$, 
of a ``construction block'' $K$, and of its boundary decomposition $\partial_{0,D} K$ and 
$\partial_{0,N} K$, we have
\[
\sigma_{DN}(\Omega_1)=\sigma_{DN}(\Omega_2)
\] 
with the account of multiplicities.
\end{thm}
\begin{proof}
The theorem is proved using the transplantation technique developed in \cite{Ber}, \cite{Bus}. 
We show that there is a one-to-one correspondence between the eigenfunctions on $\Omega_1$ and 
on $\Omega_2$. Let $u_1(x)$ be an eigenfunction of the mixed Dirichlet-Neumann boundary problem on 
$\Omega_1$. Let us represent $u_1(x)$ as follows:
\begin{equation}
u_1(x)=
\cases{
u_{11}(x), \quad & $x\in K$,\\
u_{12}(T_a x), \quad & $x\in T_a K$,
}
\end{equation}
where $u_{11}(x)$, $u_{12}(x)$ are functions on $K$ satisfying
\[
u_{11}(x)=u_{12}(x), \,\, \partial_n u_{11}(x)=-\partial_n u_{12}(x),
\] 
for $x \in \partial_a K$.
Let 
\begin{equation}
\label{lin}
u_{21}(x)=u_{11}(x)-u_{12}(x), \quad u_{22}(x)=u_{11}(x)+u_{12}(x).
\end{equation}
One can check by inspection that the function
\begin{equation}
u_2(x)=
\cases{
u_{21}(x), \quad & $x\in K$,\\
u_{22}(T_b x), \quad & $x\in T_b K$,
}
\end{equation}
is an eigenfunction of the corresponding mixed Dirichlet-Neumann boundary value problem on $\Omega_2$. It is
easy to see that inverting this procedure one obtains an eingenfunction of the problem on $\Omega_1$ from an
eigenfuntion of the corresponding problem on $\Omega_2$. Note also that since \eref{lin} is a linear
transformation we get $\sigma_{DN}(\Omega_1)=\sigma_{DN}(\Omega_2)$ with the account of multiplicities. This
completes the proof of the theorem. \end{proof}

The construction of this section also gives us the basic example of section~\ref{sec:basex} if we take the ``construction
block"  $K$ to be an isosceles right-angled triangle with legsize one, $a$ being the hypotenuse of $K$, $b$
being one of the legs,  and  $\partial_{0,N} K$ chosen to be empty.

\section{Isospectrality and multi-sheeted coverings of $K$}\label{sec:isosp} 

In this section, we indicate an alternative way of proving Theorem~\ref{thm:main}, and, at the same time, 
relate in an indirect way the spectra  $\sigma_{DN}(\Omega_j)$ and the spectra of boundary value 
problems on the ``construction block'' $K$. For illustrative purposes all the figures in this section use
a construction block $K$ with parallel sides $a$, $b$, which is different from the construction block shown 
in the previous section.

\begin{figure}[!hbt]
\begin{center}
\framebox[0.8\textwidth]{\includegraphics[width=0.60\textwidth]{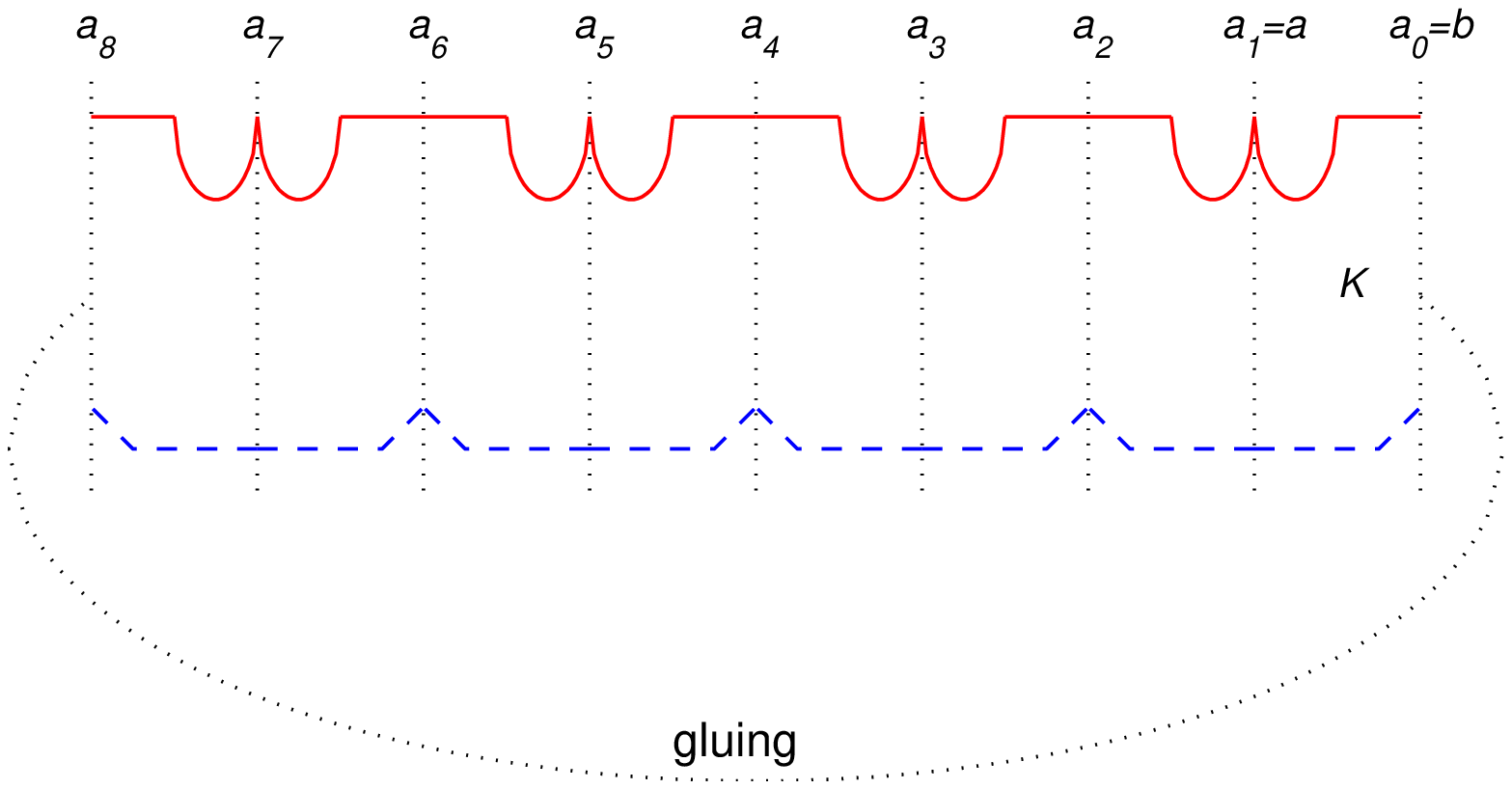}}
\caption{$\stackrel{\curvearrowright}{K}_8$, an eigth-sheeted covering of $K$. Here the construction block $K$ is bounded by 
two parallel lines $a$ and $b$. The main construction of section~\ref{sec:main} gives an isospectral pair, with $\Omega_1$ being bounded by $a_2$ 
(with the Dirichlet condition imposed) and $a_0$ (Neumann) and 
$\Omega_2$ being bounded by $a_3$ (Dirichlet) and $a_1$ (Neumann).\label{fig:cover}}
\end{center}
\end{figure}

To start with, consider an eight-sheet covering $\stackrel{\curvearrowright}{K}_8$ of the block $K$. It is constructed  
by ``gluing''
together  four copies of  $K$ and four copies of its reflection $T_a K$, and identifying the outer edges, as
shown in Figure~\ref{fig:cover}.   

The dotted lines show the original bounding lines and their images under reflections; we set 
\[
a_0:=b\,,\qquad a_1:=a\,,\qquad a_j:=T_{a_{j-1}}(a_{j-2})\quad\text{for }j=2,\dots,8
\]
and identify $a_0$ and $a_8$. 

$\stackrel{\curvearrowright}{K}_8$ is a flat manifold with boundary on which we consider a mixed Dirichlet-Neumann problem with boundary
conditions imposed according to the original choice of $\partial_{0,D} K$ and  $\partial_{0,N} K$ and their
reflections. Denote the spectrum of the corresponding Laplacian by 
$\sigma(\stackrel{\curvearrowright}{K}_8)$.

Any pair of lines $a_k\cup a_{k+4}$, $k=0,\dots,4$ defines a symmetry on $\stackrel{\curvearrowright}{K}_8$
which preserves both $a_k$ and $a_{k+4}$ and exchanges $a_{(k-j)\mod 8}$ with $a_{(k+j)\mod 8}$ for $j=1,2,3$. 
Consider, e.g., the case 
$k=0$. Then 
the lines
$a_0$, $a_4$ split $\stackrel{\curvearrowright}{K}_8$ into two identical domains; denote one of them $K_4$. The eigenfunctions 
on $\stackrel{\curvearrowright}{K}_8$ 
can be chosen in such a way that each one shall satisfy
either a Dirichlet or a Neumann boundary condition on $a_0\cup a_4$, and
\[
\sigma(\stackrel{\curvearrowright}{K}_8) = \sigma_{DD}(K_4)\cup\sigma_{NN}(K_4)\,,
\]
where $\sigma_{DD,NN}(K_4)$ stand for the spectra of the Laplacian on $K_4$ with corresponding boundary
conditions  imposed on $a_0$ and $a_4$, and the union is understood with
account of multiplicities.

Consider now these two new problems on $K_4$. For each of them, $a_2$ is the line of symmetry which divides
$K_4$ into two copies of $\Omega_1$. By the same argument,
\[
\fl\quad\sigma_{DD}(K_4) = \sigma_{DN}(\Omega_1)\cup\sigma_{DD}(\Omega_1)\qquad\text{and}\qquad
\sigma_{NN}(K_4) = \sigma_{ND}(\Omega_1)\cup\sigma_{NN}(\Omega_1)\,;
\]
here again the indices $D$ and $N$ correspond to the boundary conditions imposed on the sides $a_0$ and 
$a_2$ of $\Omega_1$.
Obviously, $\sigma_{ND}(\Omega_1)=\sigma_{DN}(\Omega_1)$ by symmetry.

Repeating the argument once more for symmetric problems on $\Omega_1$, we get
\[
\fl\quad\sigma_{DD}(\Omega_1) = \sigma_{DN}(K)\cup\sigma_{DD}(K)\qquad\text{and}\qquad
\sigma_{NN}(\Omega_1) = \sigma_{ND}(K)\cup\sigma_{NN}(K)\,,
\]
and so
\begin{equation}\label{eq:sum1}
\fl\quad\sigma_{DN}(\Omega_1)\cup\sigma_{DN}(\Omega_1)=\sigma(\stackrel{\curvearrowright}{K}_8) \setminus 
(\sigma_{DD}(K)\cup\sigma_{DN}(K)\cup\sigma_{ND}(K)\cup\sigma_{NN}(K))\,,
\end{equation}
i.e. the spectrum $\sigma_{DN}(\Omega_1)$ taken with double multiplicities is obtained by removing from the
spectrum  $\sigma(\stackrel{\curvearrowright}{K}_8)$ of the problem on the eight-sheeted covering the spectra of the four boundary value
problems on  the ``construction block'' $K$.

Let us repeat now the whole process but start with considering the symmetry with respect to $a_1\cup a_5$. Then
instead of  $K_4$ we should consider a different set $K_4'$ bounded by $a_1$ and $a_5$, which has a line of
symmetry $a_3$ dividing it into two copies of $\Omega_2$. The previous argument gives again \eref{eq:sum1},
with $\Omega_1$ replaced by $\Omega_2$, thus providing an alternative proof of  Theorem~\ref{thm:main}.

Further on, a similar argument provides us with another pair of isospectral domains with more
complicated symmetries than those of the previous section. Namely, we have
\begin{thm}
\ 
\[
\sigma_{DN}(K_4)=\sigma_{DN}(K_4')\,,
\]
or more generally
\[
\sigma_{DN}(K_{2^n})=\sigma_{DN}(K_{2^n}')\,,
\]
where $K_{2^n}$ is constructed by gluing together $2^{n-1}$ copies of $K$ and  $2^{n-1}$ copies of $T_a K$
starting with $K$, and $K_{2^n}'$ is constructed by gluing together $2^{n-1}$ copies of $K$ and  $2^{n-1}$
copies of $T_b K$ starting with $K$. 
\end{thm}

\begin{figure}[!hbt]
\begin{center}
\framebox[0.8\textwidth]{\includegraphics[width=0.60\textwidth]{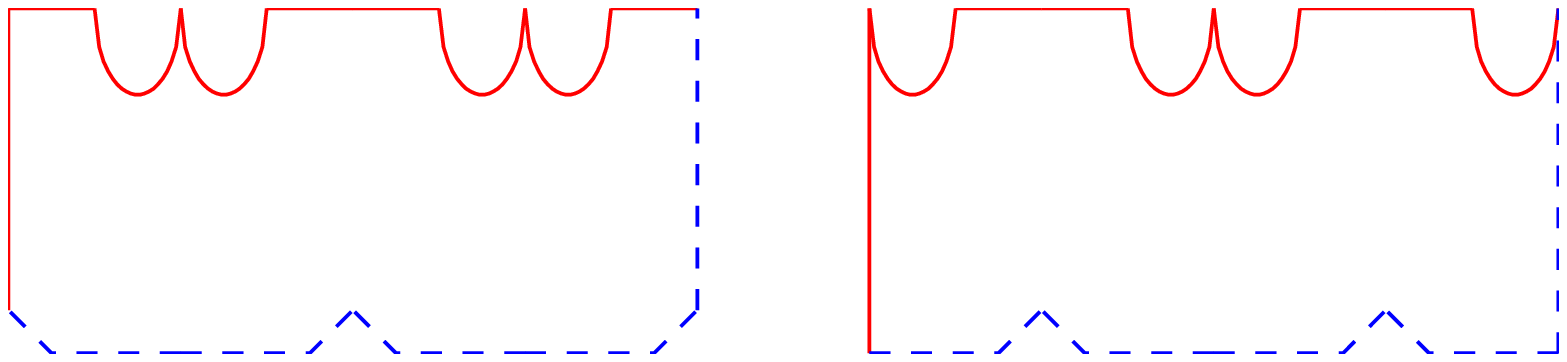}}
\caption{Isospectral problems on $K_4$ and $K_4'$. The construction block $K$ 
is the same as in Figure~\ref{fig:cover}.\label{fig:k4s}}
\end{center}
\end{figure}

See Figure~\ref{fig:k4s} for illustration in the case $n=2$. The proof is by induction in $n$ using the 
argument similar  to the one above via the construction of $2^{n+3}$-multiple covering 
$\stackrel{\curvearrowright}{K}_{2^{n+3}}$. We omit the details of the argument.

\section{Necessary conditions for mixed isospectrality}\label{sec:necess}

We present here some  necessary conditions for mixed isospectrality coming from the heat trace asymptotic 
expansion for a mixed Dirichlet-Neumann boundary problem.

\begin{thm} Let two domains $\Omega_j$, $j=1,2$, with boundary decompositions 
$\partial\Omega_j = \overline{\partial_D\Omega_j\cup\partial_N\Omega_j}$,  be isospectral in the above 
sense. 
Then the following quantities should coincide for $j=1$ and $j=2$:
\begin{itemize}
\item $\Area(\Omega_j)$;
\item $\Length(\partial_D \Omega_j)-\Length(\partial_N \Omega_j)$;
\item $\displaystyle 2\int_{\partial \Omega_j}\kappa(s)\,\mathrm{d}s
+\sum_{DD}\frac{\pi^2-\beta^2}{\beta}+
\sum_{NN}\frac{\pi^2-\beta^2}{\beta}-\frac{1}{2}\sum_{DN}\frac{\pi^2+2\beta^2}{\beta}$.
\end{itemize}
In the last formula, $\kappa$ is the curvature of the boundary, and the sums are taken over all corners of $\partial\Omega_j$ formed by 
Dirichlet-Dirichlet (DD), Neumann-Neumann (NN), and Dirichlet-Neumann (DN) parts of the boundary, 
respectively; in each case $\beta$ is a corresponding angle. 
\end{thm}
\begin{proof} The theorem follows immediately from the formulae for the first three heat trace coefficients 
in the case of a mixed Dirichlet-Neumann problem  (see \cite{Dow}) --- as the spectra coincide, the heat
trace expansions should coincide as well.
\end{proof}

\section{More examples}\label{sec:more}

In this section, we illustrate our construction by more examples. They are quite simple and do not require
explanation  apart from graphical one. In all the cases we show the construction block and resulting isospectral
domains.

\

\begin{figure}[!hbt]
\begin{center}
\framebox[0.8\textwidth]{\includegraphics[width=0.60\textwidth]{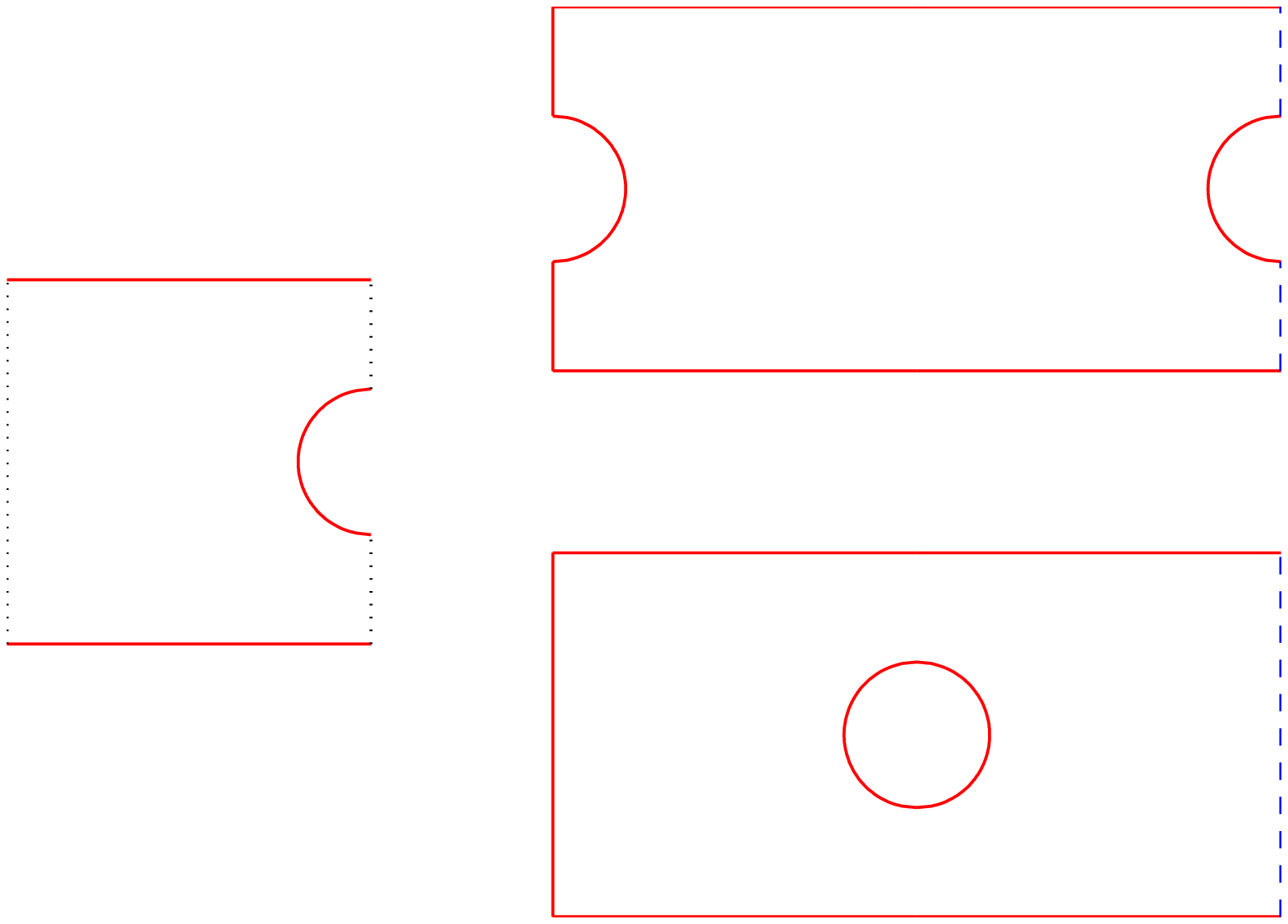}}
\caption{\small Isospectral simply connected and non-simply connected domains.\label{fig:sc_nonsc}}
\end{center}
\end{figure}

\noindent{\bf Example 1. Simply connected and not simply connected domains. } As in Figure~\ref{fig:sc_nonsc}.

\

\begin{figure}[!hbt]
\begin{center}
\framebox[0.8\textwidth]{\includegraphics[width=0.60\textwidth]{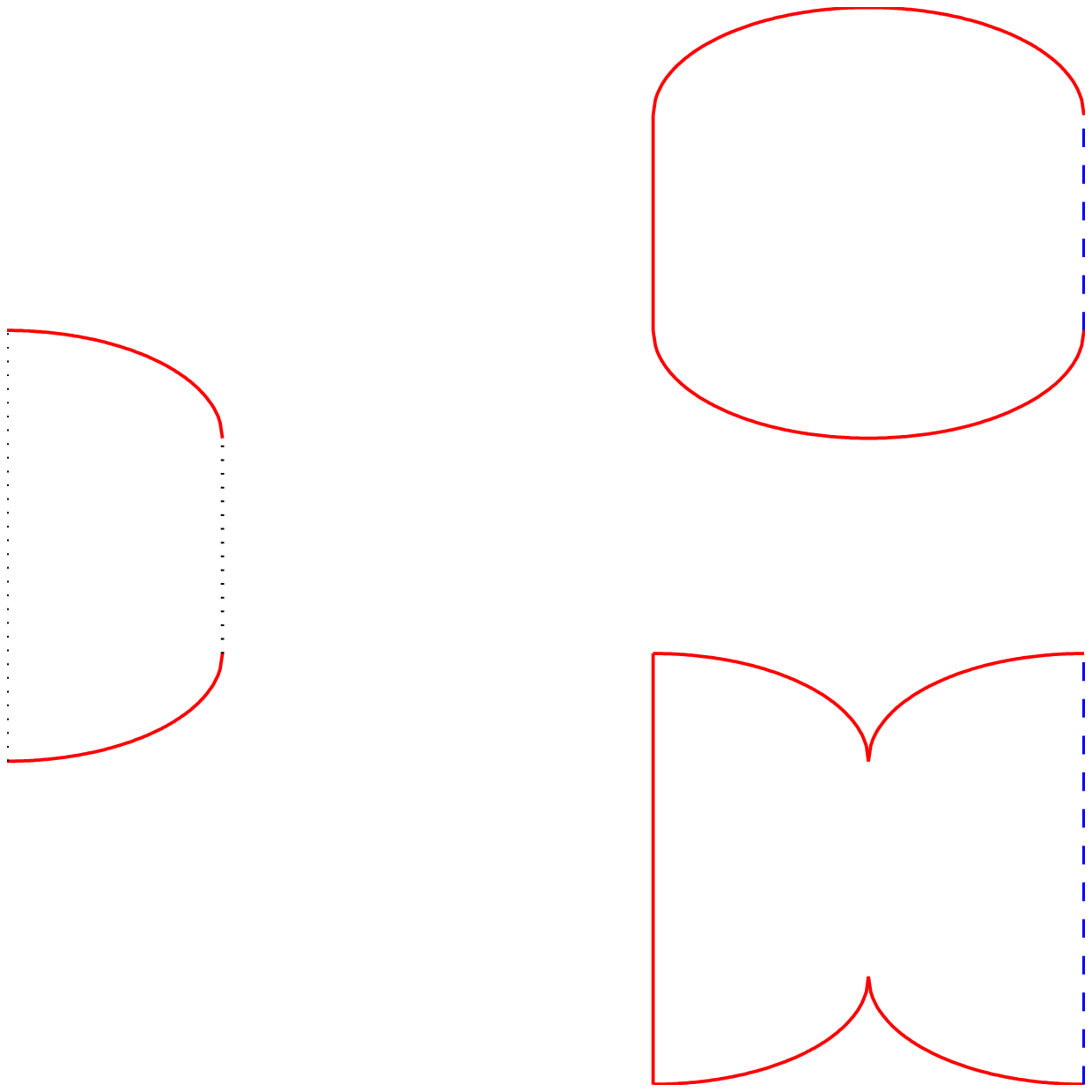}}
\caption{\small Isospectral smooth and non-smooth domains.
\label{fig:sm_nonsm}}
\end{center}
\end{figure}

\noindent{\bf Example 2. One smooth and one non-smooth domain. } As in  Figure~\ref{fig:sm_nonsm}. Note
that the boundary of the upper domain is smooth but not analytic.

\

\noindent{\bf Example 3. Four isospectral domains. } This example is more complicated than the two previous ones.

\begin{figure}[!hbt]
\begin{center}
\framebox[0.9\textwidth]{\includegraphics[width=0.90\textwidth]{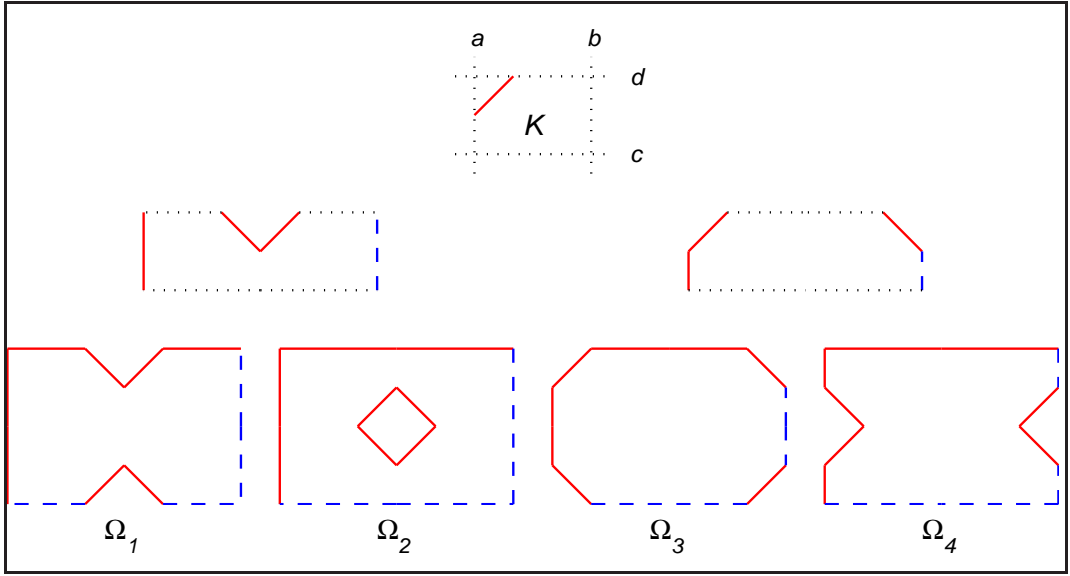}}
\caption{\small Four isospectral domains.\label{fig:fours}}
\end{center}
\end{figure}

We start with a construction block lying inside a \emph{rectangle} formed by two pairs of parallel straight lines, $a$, $b$, and $c$,
$d$, see Figure~\ref{fig:fours}. We reflect first with respect to the line $a$ and impose the Dirichlet (resp., Neumann) condition on 
$T_a b$ (resp., $b$) thus giving us the domain on the left of the middle row of Figure~\ref{fig:fours}. We now reflect it with respect to the 
lines $c$ or $d$ and impose different boundary conditions on opposite straight lines, thus giving us (by the argument of section~\ref{sec:main}) a pair of
isospectral domains $\Omega_1$ and $\Omega_2$.

If we reflect first with respect to the line $b$ and impose the Dirichlet (resp., Neumann) condition on 
$a$ (resp., $T_b a$) we obtain the domain shown on the right of the middle row of Figure~\ref{fig:fours}. 
Again reflecting with respect to $c$ and $d$, we obtain another isospectral pair $\Omega_3$ and $\Omega_4$.

It remains therefore to show that both \emph{pairs} are isospectral to each other. To do this we notice that reflecting \emph{first} with 
respect to $c$, and then with respect to $a$ and $b$ proves, by the same argument, the isospectrality of the resulting pair $\Omega_1$ and 
$\Omega_3$. Therefore, all four domains are isospectral. Note also that $\Omega_2$ is non-simply connected and the other domains are 
simply connected, thus extending Example 1.

More complicated constructions producing sets of $2^n$
isospectral non-isometric domains are possible in any dimension $n\ge 2$.

\

\begin{figure}[!hbt]
\begin{center}
\framebox[0.8\textwidth]{\includegraphics[width=0.60\textwidth]{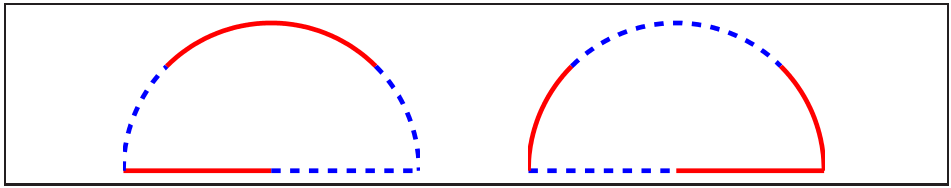}}
\caption{\small Two problems on a half-disk.\label{fig:halfs}}
\end{center}
\end{figure}

\noindent{\bf Example 4. Domains isospectral with respect to Dirichlet-Neumann swap.} This is discussed at
length in \cite{JLNP}, where we construct a series of examples of boundary value problems whose spectrum
remains the same when swapping  Dirichlet and Neumann parts of the boundary. The main example is of two
problems on a half-disk, see Figure~\ref{fig:halfs}.

This example has an unexpected application \cite{JLNNP}, which has in fact first attracted our attention to
Dirichlet-Neumann isospectrality. 

\section{Can one count the shape of a drum?}\label{sec:see}

Very recently, Uzy Smilansky and collaborators have noticed \cite{GnSmSo, ShSm} that in some situations pairs of previously known 
isospectral objects (manifolds or graphs) can be distinguished by comparing their \emph{nodal sequences}.  Namely, if $\lambda_m$ is the $m$-th eigenvalue
and $u_m$ is a corresponding eigenfunction of the Dirichlet Laplacian on $\Omega$, denote by $N_m(\Omega):=\{x\in\Omega: u_m(x)=0\}$ 
the nodal set of $u_m$ and by $\nu_m(\Omega)$ the number of connected components of $\Omega\setminus N_m(\Omega)$. Note that $\nu_m(\Omega)$
is only well defined in such a way for a \emph{simple} eigenvalue $\lambda_m$ and that the case of multiple eigenvalues 
requires some special treatment.

Numerical experiments in \cite{GnSmSo, ShSm} based on statistical analysis of a normalized nodal domain count, see \cite{BlGnSm}, showed that 
isospectral pairs provide in fact different distributions of nodal domains.

\begin{figure}[!hbt]
\begin{center}
\framebox[0.8\textwidth]{\includegraphics[width=0.60\textwidth]{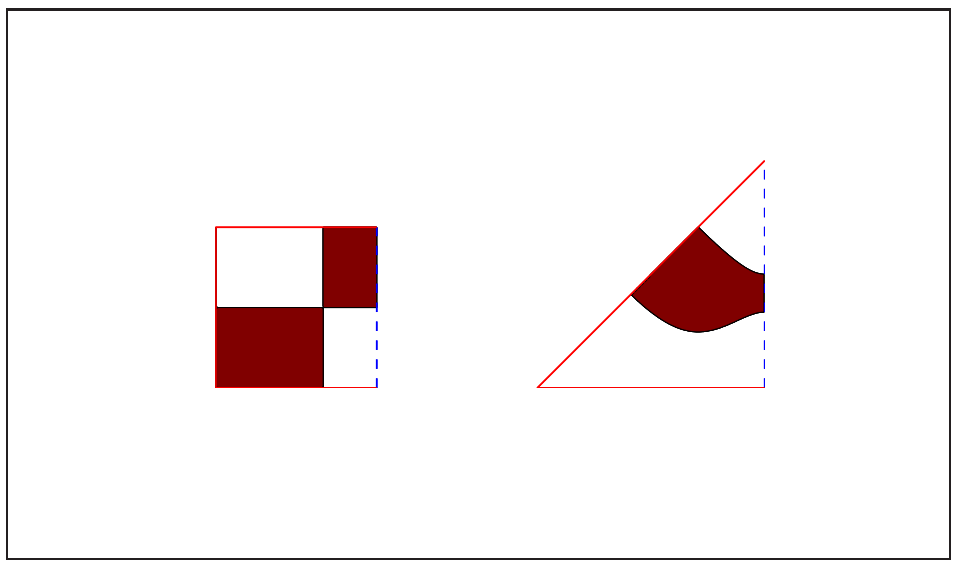}}
\caption{\small Nodal domains for the basic example of section~\ref{sec:basex}.\label{fig:levels}}
\end{center}
\end{figure}

We note here that, at least for some low eigenvalues, the pairs of isospectral domains with mixed boundary conditions 
considered in this paper also produce different nodal count. Consider, e.g. our basic example of a square and a triangle in section~\ref{sec:basex}.
The nodal domains for the eigenfunctions corresponding to the fourth eigenvalues $\lambda_{1,2}=\mu_{3,0}$ are shown in Figure~\ref{fig:levels},
and we see that $\nu_4(\Omega_1)=4\ne 3=\nu_4(\Omega_2)$.

As finite element calculations show, the nodal domain count is also different, e.g., for the third eigenfunctions of the two problems in the half-disk isospectral with respect to the
Dirichlet-Neumann swap (Example 4 of the previous section) with four and two nodal domains, respectively.

It would be of course interesting to prove that \emph{any} two isospectral mixed Dirichlet--Neumann problems can be distinguished by their nodal 
count; this question remains open.



\ack
The authors are grateful to Uzy Smilansky for valuable discussions. 
The research of I.P. was partially supported by NSERC and  FQRNT.

\end{document}